\documentclass[a4paper, 12pt, titlepage, roman]{article}
\usepackage{MRHStyle}
\usepackage{JournalOfExpMath}

\usepackage[square]{natbib}
\setcitestyle{citesep={;}}

\begin{document}
\title{\Huge \textbf Collatz polynomials: an introduction with bounds on their zeros}
\author{\Large Matt Hohertz \\ \normalsize Corresponding author \\ Department of Mathematics, Rutgers University \\ mrh163@math.rutgers.edu \\ \\ \Large Bahman Kalantari \\ \normalsize Co-author \\ Department of Computer Science, Rutgers University \\ kalantari@cs.rutgers.edu \\ \\}

\maketitle

\begin{abstract}
The Collatz Conjecture (also known as the 3x+1 Problem) proposes that the following algorithm will, after a certain number of iterations, always yield the number 1: given a natural number, multiply by three and add one if the number is odd, halve the resulting number, then repeat.  In this article, for each $N$ for which the Collatz Conjecture holds we define the $N^{th}$ Collatz polynomial to be the monic polynomial with constant term $N$ and $k^{th}$ term (for $k > 1$) the $k^{th}$ iterate of $N$ under the Collatz function.  In particular, we bound the moduli of the roots of these polynomials, prove theorems on when they have rational integer roots, and suggest further applications and avenues of research.
\end{abstract}

\bd{Keywords} \quad Collatz conjecture, generating function, 3x+1, polynomial roots, bounds on roots \par
\bd{Word count} \quad 1887

\section{Introduction}
The Collatz Conjecture (or 3x+1 Problem) proposes that repeated iteration of any natural number $N$ under the \ic{Collatz function} 
\begin{equation}
c(N) =
\begin{cases}
	\frac{3N+1}{2}, & N\mbox{ odd} \\
	\frac{N}{2}, & N\mbox{ even}
\end{cases}
\end{equation} 
or, equivalently, under the commonly-used \ic{alternative Collatz function}
\begin{equation}
\widehat{c}(N) = 
\begin{cases}
	3N+1, & N\mbox{ odd} \\
	\frac{N}{2}, & N\mbox{ even}
\end{cases}
\end{equation} 
eventually leads to the number 1.   Though easily stated and understood, this conjecture has proved notoriously intractable, to an extent that even Paul Erdős famously remarked of it that ``mathematics is not yet ready for such problems" \citep{lag1}. The failure of direct approaches to the problem, as well as the uselessness of experimental approaches for providing its rigorous solution, has led to its rephrasing as a problem in other domains, including ergodic theory \citep{matthews} and the theory of computation \citep{conway}, as well as the solution of more approachable related problems (as exemplified by \citep{tao}).  \par
	To each natural number $N$ for which the Collatz Conjecture holds, we define the $N^{th}$ \ic{Collatz polynomial} to be that monic polynomial whose coefficients are iterates of the Collatz function.  We then use bounds by Kalantari and showcased in \citep{kalantari_2004} and \citep{kalantari}, as well as Fujiwara's bound, to bound the moduli of the roots of such polynomials. The upper bound we find improves asymptotically with the base $t$ of $N$ (a parameter which we will define) and the lower bound improves asymptotically with $N$ itself.

\section{Definitions}\label{sec:defs}
Define the \emp{$N^{th}$ Collatz polynomial} to be the generating polynomial 
\begin{align}
	P_N(z) &:= \sum_{k\geq 0} c^k(N)z^k \\
	&= N + c(N)z + ... + c^j(N)z^j + ... + 2^tz^{n - t} + ... + 2z^{n-1} + z^n
\end{align} 
of iterations of $N$ under the Collatz function, where
\begin{itemize}
	\item $c^k:\N\rightarrow\N$ is defined recursively to be
	\begin{equation}
	c^k(N) = 
	\begin{cases}
		N, & k = 0 \\
		c(c^{k-1}(N)), & 1 \leq k \leq n \\
		0, & k > n \mbox{ (where }n\mbox{ is defined as below)} \\
	\end{cases}
	\end{equation}
	\item $n$ is the \emp{total stopping time} of $N$, the least $k$ such that 
	\begin{equation}
		c^k(N) = 1
	\end{equation}
	\item $t$ is the \emp{base} of $N$, the logarithm base 2 of the first iterate of $N$ that is a power of two.  
\end{itemize}
\paragraph{\textsc{Example}} $P_5(z) = 5 + 8z + 4z^2 + 2z^3 + z^4.$ \par
In section \ref{sec:altcpol}, we will also consider the $N^{th}$ \ic{alternative Collatz polynomial} $\widehat{P}_N(z)$, whose coefficients are the iterates of $N$ under the alternative Collatz function.
\paragraph{\textsc{Example}} $\widehat{P}_5(z) = 5 + 16z + 8z^2 + 4z^3 + 2z^4 + z^5.$
%\paragraph{\textsc{remark}} Note that, by considering only monic polynomials, we are assuming the Collatz conjecture to be true.

\section{Bounds for the moduli of roots}
	\subsection{An upper bound}
	
		\begin{lemma} \label{lemma:unisol} 
		For $m\geq 2$, the equation
		\begin{equation}
			f(x) = x^{m-1}+x-1 = 0 \label{eq:mpol}
		\end{equation}
		has exactly one solution $r_m$ on the interval $[\frac{1}{2},1)$.
		\begin{proof}
		The derivative
		\begin{equation}
			f'(x) = \inps{m-1}x^{m-2} + 1
		\end{equation}
		of the left-hand side is positive on the interval, so the solution to the equation is unique if it exists.  If $m = 2$ then the solution is $x = \frac{1}{2}$; otherwise,
		\begin{equation}
			f\inps{\frac{1}{2}} = \inps{\frac{1}{2}}^{m-1} + \frac{1}{2} - 1 < 0
		\end{equation}
		\begin{equation}
			f(1) = 1 > 0
		\end{equation}
		and the Intermediate Value Theorem implies the existence of a solution.
		\end{proof}
		\end{lemma}
		
		\begin{lemma} \label{lemma:powbnd}
		\begin{equation*}
			c^{j}(N)\leq 2^{n-j}
		\end{equation*}
		\begin{proof}
		The inequality holds for $j = n$, and $c^{j-1}(N)$ is either exactly twice $c^{j}\inps{N}$ or strictly less than $c^{j}\inps{N}$.
		\end{proof}
		\end{lemma}
		
		\begin{lemma} 
		If $t\geq 3$ then \label{lemma:fracbnd}
			$$\frac{2^{t+2} + 1}{3} \leq \frac{11}{64}\cdot 2^{t+3}$$
		\begin{proof}
			Equality holds for $t = 3$, and from the $t^{th}$ to the $\inps{t+1}^{th}$ term the left-hand side increases by a ratio of strictly less than 2.
		\end{proof}
		\end{lemma}
		
		\begin{lemma}
		The quantity $r_m$ increases to 1 as $m$ increases.  \label{lemma:rmdec}
		\begin{proof}
		We rewrite the equation
		$$x^{m - 1} + x - 1 = 0$$
		as 
		\begin{equation}
			x = 1 - x^s
		\end{equation}
		and take the derivative of $x$ with respect to $s$.
		\begin{equation}
			\frac{dx}{ds} = -\inps{\log x}x^s
		\end{equation}
		By Lemma \ref{lemma:unisol}, $\log x < 0$, and so the derivative is positive and $r_m$ is increasing. \par
		To prove that $r_m\rightarrow 1$, we first rearrange Equation \eqref{eq:mpol} and take logarithms to obtain
		\begin{align}
			\log\inps{r_m} &= \frac{\log\inps{1 - r_m}}{m - 1}
		\end{align}
		By the Monotone Convergence Theorem, $r_m$ approaches a finite limit less than or equal to 1 as $m\rightarrow\infty$.  Suppose $r_m\not\rightarrow 1$.  Then $\frac{\log\inps{1 - r_m}}{m - 1}\rightarrow 0$ as $m\rightarrow\infty$.  However, $\log\inps{r_m}\not\rightarrow 0$ as $m\rightarrow\infty$, a contradiction.  Therefore, $r_m\rightarrow 1$.
		\end{proof}
		\end{lemma}
		
		\begin{theorem}\label{thm:upperbnd}
		Given $N\in\N$ with base $t$, for any root $\xi$ of $P_N$, 
		$$\boxed{|\xi|\leq h(t):=\frac{2}{r_{t+3}}\cdot\inps{\frac{75}{32} + 2^{-(t+2)}}^{\frac{1}{t+2}}}$$
		where $r_{t+3}$ is the unique root of
		\begin{equation}
			x^{t+2}+x-1
		\end{equation}
		contained in the interval $x\in[\frac{1}{2},1)$. \par Moreover, $h(t)$ decreases monotonically to 2 as $t\rightarrow \infty$.
		\begin{proof}
		Let
		\begin{equation}\label{eq:genpol}
			f(z) = a_nz^n + a_{n-1}z^{n-1} + ... + a_1z + a_0
		\end{equation}
		be a general complex polynomial of degree $n$ of which $\xi$ is any root; for any natural number $m\geq 2$, define
		\begin{equation}
			A_{m,k} = 
			\begin{pmatrix}
			a_{n-1} & a_{n-2} & \cdots & a_{n-m+1} & a_{n-k+1} \\
			a_{n} & a_{n-1} & \cdots & a_{n-m+2} & a_{n-k+2} \\
			0 & a_{n} & \cdots & a_{n-m+3} & a_{n-k+3} \\
			\vdots & \ddots & \ddots & \vdots & \vdots \\
			0 & \cdots & 0 & a_n & a_{n-k+m} \\
			\end{pmatrix}
		\end{equation}
		Kalantari proves in \citep{kalantari_2004} and \citep{kalantari} that any root $\xi$ satisfies
		\begin{equation}\label{eq:polbnd}
			|\xi|\leq\frac{1}{r_m}\cdot\max_{k: \linebreak m\leq k \leq m + n - 1}\aval{\det\inps{A_{m-1,k}}}^{1/(k-1)} =: U_m
		\end{equation} 
		where $r_m$ is the root of 
			$$x^{m-1} + x - 1$$
		contained in the interval $\left[\frac{1}{2},1\right)$ (a root that, by Lemma \ref{lemma:unisol}, exists and is unique).  Since $N$ has base $t$, the preimage of $2^t$ under $c$ is $\frac{2^{t+1}-1}{3}$.  Thus, setting $m = t+3$ we obtain, for the Collatz polynomial
		$$P_N(z) = \sum_{j=0}^n a_jz^j$$ 
		the matrix
		\begin{equation} 
			A_{t+2,k} = \begin{pmatrix}
			2 & 4 & \cdots & 2^t & \frac{2^{t+1} - 1}{3} & a_{n-k+1} \\
			1 & 2 & \cdots & 2^{t-1} & 2^t  & a_{n-k+2} \\
			0 & 1 & 2 & \cdots & 2^{t-1} & a_{n-k+3} \\
			\vdots & \ddots & \ddots & \ddots & \vdots & \vdots \\
			\vdots & & \ddots & 1 & 2 & a_{n-k+t+2}\\
			0 & \cdots & \cdots & 0 & 1 & a_{n-k+t+3}\\
			\end{pmatrix}
		\end{equation}	
		Subtracting twice the second row from the first row yields a matrix with the same determinant; thus,
		\begin{align} 
			\aval{\det\inps{A_{t+2,k}}}&= \aval{\det{
			\begin{pmatrix}
			0 & 0 & \cdots & 0 & -\frac{2^{t+2} + 1}{3} & a_{n-k+1} - 2a_{n-k+2} \\
			1 & 2 & \cdots & 2^{t-1} & 2^t  & a_{n-k+2} \\
			0 & 1 & 2 & \cdots & 2^{t-1} & a_{n-k+3} \\
			\vdots & \ddots & \ddots & \ddots & \vdots & \vdots \\
			\vdots & & \ddots & 1 & 2 & a_{n-k+t+2}\\
			0 & \cdots & \cdots & 0 & 1 & a_{n-k+t+3}\\
			\end{pmatrix}
			}} \\
			&= \aval{\inps{-1}^{t+2}\inps{ -\frac{2^{t+2} + 1}{3}\cdot a_{n-k+t+3} + 2a_{n-k+2}-a_{n-k+1}}} \label{eq:sumofcoeffs} \\
			&= 
			\begin{cases}
				\frac{2^{t+2} + 1}{3}\cdot a_{n-k+t+3}, &a_{n-k+1}\mbox{ even and non-zero} \\
				\frac{2^{t+2} + 1}{3}\cdot a_{n-k+t+3}+2a_{n-k+2}, &a_{n-k+1}\mbox{ zero} \\
				\frac{2^{t+2} + 1}{3}\cdot a_{n-k+t+3}+2a_{n-k+1} + 1, &a_{n-k+1}\mbox{ odd} \\
			\end{cases} \label{eq:avcases}
		\end{align} 
		By Lemmas \ref{lemma:powbnd} and \ref{lemma:fracbnd}, these bounds become, respectively,
		\begin{equation}
			\leq\begin{cases}
			\frac{11}{32}\cdot 2^{k-1}, &a_{n-k+1}\mbox{ even and non-zero} \\
			\frac{43}{32}\cdot 2^{k-1}, &a_{n-k+1}\mbox{ zero} \\								
			\inps{\frac{75}{32} + 2^{-(k-1)}}\cdot 2^{k-1}&a_{n-k+1}\mbox{ odd} \\
			\end{cases}
		\end{equation}
		Hence the quantity $U_{t+2}$ of Equation \eqref{eq:polbnd} is bounded above by 
		$$\frac{2}{r_{t+3}}\cdot\inps{\frac{75}{32} + 2^{-(k-1)}}^{1/(k-1)}$$
		For $k > 1$, the derivative of this quantity with respect to $k$ is negative, yielding the upper bound
		\begin{equation} \label{eq:boundint}					
			\frac{2}{r_{t+3}}\cdot\inps{\frac{75}{32} + 2^{-(t+2)}}^{1/(t+2)}
		\end{equation}
		
		Since $\inps{\frac{75}{32} + 2^{-(t+2)}}^{1/(t+2)}\rightarrow 1$ as $t\rightarrow\infty$ and has negative derivative with respect to $t$, it follows that $h(t)$ decreases to 2 as $t\rightarrow\infty$. 
		\end{proof}
		\end{theorem}
		Values of $h(t)$ for increasing values of $t$ are shown below: \par
\begin{table}[h]
	\begin{tabular}{|c|c|c|c|c|c|}
		\hline
		$t$ & 3 & 10 & $10^3$ & $10^5$ & $10^7$ \\ \hline
		$h(t)$ & 3.1498 & 2.5185 & 2.0122 & 2.0002 & 2.000002\\
		\hline
	\end{tabular}
	\caption{Particular values of $h(t)$}
\end{table}
		
	\subsection{A lower bound}
		\begin{lemma} \label{lemma:ratbnd}
			$$\frac{c^j(N)}{N} \leq \inps{\frac{3}{2}}^j\cdot\inps{1 + \frac{1}{N}}$$
			\begin{proof}
				A direct corollary of Lemma 5 in \citep{berg}.
			\end{proof}
		\end{lemma}
		\begin{lemma}[Fujiwara's bound, \citep{fujiwara}]
			For any root $\xi$ of the general complex polynomial of Equation \eqref{eq:genpol},
			\begin{equation*}
				\aval{\xi}\leq \max\Bigg\{\max_{i=1,...,n-1}\aval{\frac{a_{n-i}}{a_n}}^{\frac{1}{i}}, \aval{\frac{a_0}{2a_n}}\Bigg\}
			\end{equation*}
		\end{lemma}
		\begin{theorem} \label{thm:lb}
		For any root $\xi$ of $P_N$, 
			$$\boxed{|\xi| \geq \frac{2}{3\inps{1 + \frac{1}{N}}}}$$
			\begin{proof}
				We apply Fujiwara's bound to the reciprocal polynomial of $P_N(z)$.  By Lemma \ref{lemma:ratbnd}, 
				\begin{equation}
					\aval{\frac{a_i}{a_0}}^{\frac{1}{i}} \leq \frac{3}{2}\cdot \inps{1 + \frac{1}{N}}^{\frac{1}{i}}
				\end{equation}
				The maximum of this quantity over $i = 1, ..., n-1$ is $\frac{3}{2}\inps{1+\frac{1}{N}}$, which is greater than 
				\begin{equation}
					\aval{\frac{a_n}{2a_0}} = \frac{1}{2N}
				\end{equation}
				Since the zeros of the original polynomial are the reciprocals of those of the reciprocal polynomial, the bound follows.
			\end{proof}
		\end{theorem}
			
\section{Remarks}\label{sec:compar}
\begin{enumerate}
	\item The upper bound given by Theorem \ref{thm:upperbnd} could be expressed in terms of $N$ (albeit perhaps as a much less tight bound) if there existed a function $\ell\inps{N}$ bounding $t$ below uniformly -- \ic{i.e.}, such that
		\begin{equation}
			\ell\inps{N}\leq t
		\end{equation}
		However, any such function must have the property
		\begin{equation}
			\liminf_{N\rightarrow} \ell(N) = 3
		\end{equation}
		because, in particular, $2^k\cdot 5$ has base 3 for all $k$.
	\item \label{item:optimbnd} As proved in \citep{jin}, the bound from Equation \eqref{eq:polbnd} converges to the greatest absolute value of the roots of the polynomial $f(z)$ as $m\rightarrow\infty$.  Thus, we would expect that taking a greater number of iterates of $N$ into account would result in a more precise bound for $P_N(z)$.  However, such bounds become difficult to calculate by hand.  In calculating $|\det(A_{m-1,k})|$ for $m = t + 2 + \ell$, there are $2^{\ell-1}$ possible preimages of $2^t$ under $c$, hence $2^{\ell-1}$ different determinants to calculate and compare.  Moreover, the technique exhibited in the proof of Theorem \ref{thm:upperbnd} does not extend to the calculations of lower bounds: whereas 
	\begin{equation}
	\frac{a_{n-1}}{a_n} = \frac{a_{n-2}}{a_{n-1}} = \cdots = \frac{a_{n-t}}{a_{n-t+1}} = 2
	\end{equation} 
	there is no common ratio between the consecutive terms $a_0, a_1, ...$. 
	\item Applied to the reciprocal polynomial of $P_N(z)$, $U_m$ bounds $|\xi|$ below poorly; applied to $P_N(z)$ itself, however, it bounds $|\xi|$ above better than other standard bounds, even at the relatively low values of $m$ we use.\par
For example, Fujiwara's bound provides at best the unhelpful inequality $|\xi| < 2\cdot\inps{\frac{3}{2}}^{n-1}\inps{N+1}$.  Indeed, 
\begin{equation}
	\aval{\frac{a_0}{2a_n}} = \frac{N}{2}
\end{equation}
while Lemma \ref{lemma:powbnd} implies
\begin{equation}
	\aval{\frac{a_{n-i}}{a_n}}^{1/i} \leq 2
\end{equation}
and Lemma \ref{lemma:ratbnd} implies
\begin{equation}
	\aval{\frac{a_{n-i}}{a_n}}^{1/i} \leq \inps{\frac{3}{2}}^{\frac{n}{i}-1}\inps{N+1}^{1/i} < \inps{\frac{3}{2}}^{n-1}\inps{N + 1}
\end{equation}
	\end{enumerate}

%\label{lemma:powbnd} \par
Moreover, the following bound by Sun and Hsieh is superior to Cauchy's bound \citep{sun_hsieh} and improves on the bound given by Fujiwara, yet is strictly weaker than ours for any value of $t$.
\begin{lemma}[Sun-Hsieh bound, \citep{sun_hsieh}]
Let
		\begin{equation}
			f(z) = a_nz^n + a_{n-1}z^{n-1} + ... + a_1z + a_0 
		\end{equation}
be a general complex polynomial of degree $n$ of which $\xi$ is any root.  Then
\begin{equation}
	|\xi|\leq 1 + d
\end{equation}
where $d$ is the unique positive real root of
\begin{equation}
	x^3 + \inps{2-\aval{a_{n-1}}}x^2 + \inps{1-\aval{a_{n-1}}-\aval{a_{n-2}}}x - \max\sbld{\aval{a_i}}
\end{equation}
\end{lemma}

Since $a_{n-1} = 2$ and $a_{n-2} = 4$ for $N > 2$, we have that $d$ is the positive real root of
\begin{equation}\label{eq:cubic}
	x^3-5x- \max\sbld{\aval{a_i}} = 0
\end{equation}
\begin{proposition}
	Let $d(N)$ signify the value of $d$ yielded by $P_N(z)$.  Then $$1 + d(N) > h(3) = 3.1498$$
	\begin{proof}
		Let $d_1(t)$ signify the positive real root of
		\begin{equation}\label{eq:cubict}
			x^3 - 5x - t = 0
		\end{equation}
		  Then $d_1'(t) > 0$.  Indeed, we rearrange Equation \eqref{eq:cubict} and take the derivative with respect to both sides to obtain
		\begin{align}
			\frac{dx}{dt} &= \frac{1}{5}\inps{3x^2\cdot \frac{dx}{dt} - 1} \\
			\frac{dx}{dt} &= \frac{1}{3x^2-5}
		\end{align}
		Since $d_1(t) > \sqrt{\frac{5}{3}}$ for all $t > 0$, this quantity is positive, hence 
		\begin{equation}
			d(N) = d_1\inps{\max_k c^k(N)} \geq d_1(N) > d_1(0)
		\end{equation}
		Since $1 + d_1(0) = 1 + \sqrt{5} = 3.23606\cdots$, we are done.
	\end{proof}
\end{proposition}

\section{Bounds for the alternative Collatz polynomial} \label{sec:altcpol}
	Define the \ic{alternative Collatz polynomial} $\widehat{P}_N$ to be the generating polynomial whose coefficients are iterates of the alternative Collatz function, often itself designated the ``Collatz function,"
	\begin{equation}
	\widehat{c}(N) = \begin{cases}
	3N+1, & N\mbox{ odd} \\
	\frac{N}{2}, & N\mbox{ even}
	\end{cases}
	\end{equation} 
	\begin{lemma}
		$$\frac{1}{3}\inps{5\cdot 2^t + 1} \leq \frac{41}{96}\cdot 2^{t+2}$$
		
		\begin{proof}
			This proof is an analogue to the proof of Lemma \ref{lemma:fracbnd}.  Equality holds for $t = 3$, and the left-hand side increases by a ratio of less than 2 as $t$ increases by 1.
		\end{proof}
	\end{lemma}
	\begin{theorem}
		For any root $\xi$ of $\widehat{P}_N$,
		$$|\xi| \leq \frac{2}{r_6}\cdot \inps{\frac{521}{96} + 2\cdot 2^{-(t+2)}}^{\frac{1}{t+2}}\leq 3.72444268658138218$$
		\begin{proof}
			Analogous to the proof of the upper bound of $P_N$.  The quantity in the absolute values of Equation \ref{eq:sumofcoeffs} becomes
			\begin{equation}
			-\frac{1}{3}\inps{5\cdot 2^t + 1}a_{n-k+t+3} + 2a_{n-k+2} - a_{n-k+1}
			\end{equation}
			and the cases of Equation \eqref{eq:avcases} become
			\begin{equation}
			\begin{cases}
			\frac{1}{3}\inps{5\cdot 2^t + 1}\cdot a_{n-k+t+3}, &a_{n-k+1}\mbox{ even and non-zero} \\
			\frac{1}{3}\inps{5\cdot 2^t + 1}\cdot a_{n-k+t+3}+2a_{n-k+2}, &a_{n-k+1}\mbox{ zero} \\
			\frac{1}{3}\inps{5\cdot 2^t + 1}\cdot a_{n-k+t+3}+5a_{n-k+1} + 2, &a_{n-k+1}\mbox{ odd} \\
			\end{cases}
			\end{equation}
			These quantities are bounded above by, respectively,
			\begin{equation}
			\leq\begin{cases}
			\frac{41}{96}\cdot 2^{k-1}, &a_{n-k+1}\mbox{ even and non-zero} \\
			\frac{137}{96}\cdot 2^{k-1}, &a_{n-k+1}\mbox{ zero} \\
			\inps{\frac{521}{96} + 2\cdot 2^{-(k-1)}}\cdot 2^{k-1}, &a_{n-k+1}\mbox{ odd}
			\end{cases}
			\end{equation}
		\end{proof}
	\end{theorem}
	\begin{lemma}\label{lemma:ratbndalt}
		$$\frac{\widehat{c}^j(N)}{N}\leq 3^j\cdot\inps{1 + \frac{1}{N}}$$
		\begin{proof}
			Analogous to the proof of Lemma 5 of \citep{berg}.  Since
			\begin{equation}
				\widehat{c}^{m+1}(N) \leq 3\cdot \widehat{c}^m(N) + 1
			\end{equation}
			it follows that 
			\begin{equation}
				\widehat{c}^{m+1}(N) + 1\leq 3\cdot\inps{\widehat{c}^m(N) + 1}
			\end{equation}
			From the initial condition $\widehat{c}^0(N) = N$, the conclusion follows.
		\end{proof}
	\end{lemma}
	\begin{theorem}
		For any root $\xi$ of $\widehat{P}_N$,
		$$|\xi|\geq \frac{1}{3\inps{1+\frac{1}{N}}}$$
		\begin{proof}
			Fujiwara's bound combined with Lemma \ref{lemma:ratbndalt}.
		\end{proof}
	\end{theorem}

\section{Miscellaneous theorems}

	\begin{theorem}
		For $N\geq 3$, $P_N$ has at least one non-real root.
	
		\begin{proof}
			First, note that $P_N$ has no non-negative roots: this fact is a consequence of Descartes's Rule of Signs and the fact that $P_N$ has a non-zero constant term.  Suppose that all roots of $P_N$ are real and label them in increasing order:
			\begin{equation}
				-2\leq r_1 \leq r_2 \leq ... \leq r_n < 0
			\end{equation}
			where $-2\leq r_1$ follows from the fact that all the roots are negative and have sum $-2$.  By Vieta's formulas,
			\begin{equation}
				3 \leq N = \prod |r_i| 
			\end{equation}
			If $-1\leq r_1$, then 
			\begin{equation}
				3 \leq N = \prod |r_i| \leq 1
			\end{equation}
			a contradiction.  And if $-2 \leq r_1 < -1$, then $-1 < r_2 \leq ... \leq r_n < 0$, yielding
			\begin{equation*}
				\frac{3}{2} \leq \frac{N}{2} \leq \frac{N}{|r_1|} = \prod_{i=2}^n |r_i| < 1
			\end{equation*}
			also a contradiction.
	\end{proof}
\end{theorem}

Hereafter, let $m(N)$ (just $m$ when context allows) be the number of odd numbers in the Collatz trajectory of $N$, \ic{i.e.}, the sequence
$$N, c(N), ..., c^{n}(N) = 1$$

\begin{lemma}
	If $N_i$ is the $i^{th}$ odd number appearing in the Collatz trajectory of $N$, with each $N_i$ appearing at the end of a subsequence
	$$2^{\ell_i-1}N_i, 2^{\ell_i-2}N_i, ..., N_i$$
	of length $\ell_i$, then
	$$P_N(z) = \sum_{k = 1}^{m(N)} \inps{2^{\ell_k}N_k}\inps{z^{\ell_1 + ... +\ell_{k-1}}}\inps{\frac{1-\inps{\frac{z}{2}}^{\ell_k}}{2-z}}$$
\end{lemma}
\iffalse
\begin{proof}
	If $m = 1$ then the generating function is
	$$2^{\ell-1}N_1 + 2^{\ell-2}N_1z + ... + N_1z^{\ell-1} = 2^{\ell_1 - 1}N_1\cdot\frac{1-\inps{\frac{z}{2}}^{\ell_1}}{1 - \frac{z}{2}}$$
	Suppose that formula holds for $m$, and that $N$ has $m + 1$ odd numbers in its sequence.  Then $c(N_1)=2^{\ell_2}N_2$, a quantity with $m$ odd numbers in its sequence.  Moreover,
	\begin{align}
	P_N(z) &= 2^{\ell_1-1}N_1 + 2^{\ell_1-2}N_1z + ... + N_1z^{\ell_1-1} + z^{\ell_1}\cdot P_{2^{\ell_2}N_2}\inps{z} \\
	&= 2^{\ell_1 - 1}N_1\cdot\frac{1-\inps{\frac{z}{2}}^{\ell_1}}{1-\frac{z}{2}} + z^{\ell_1}\inps{ \sum_{k = 2}^{m+1} 2^{\ell_k}N_k \cdot z^{\ell_2 + ... \ell_{k-1}}\cdot \frac{1-\inps{\frac{z}{2}}^{\ell_k}}{2 - z}} \\
	&= 2^{\ell_1}N_1\cdot\frac{1-\inps{\frac{z}{2}}^{\ell_1}}{2-z} + \inps{ \sum_{k = 2}^{m+1} 2^{\ell_k}N_k \cdot z^{\ell_1 + ... \ell_{k-1}}\cdot \frac{1-\inps{\frac{z}{2}}^{\ell_k}}{2-z}} \\
	&= 2^{\ell_1}N_1\cdot\frac{1-\inps{\frac{z}{2}}^{\ell_1}}{2-z} + \inps{ \sum_{k = 2}^{m+1} 2^{\ell_k}N_k \cdot z^{\ell_1 + ... \ell_{k-1}}\cdot \frac{1-\inps{\frac{z}{2}}^{\ell_k}}{2 - z}}	\\	
	&=  \sum_{k = 1}^{m+1} 2^{\ell_k}N_k \cdot z^{\ell_1 + ... \ell_{k-1}}\cdot \frac{1-\inps{\frac{z}{2}}^{\ell_k}}{2 - z}
	\end{align}
\end{proof} \fi
\begin{theorem}
	With notation as above, $P_N(-2) = 0$ if and only if $\ell_i$ is even for all $i=1,...,m$.
	
	\begin{proof}
		One direction follows by direct substitution of $-2$ into the formula for $P_N(z)$.  \par
		For the other direction, note that the equality
		\begin{align}
		0 = P_N(-2) &= \sum_{k=1}^m \inps{2^{\ell_k}N_k}\inps{(-2)^{\ell_1+...+\ell_{k-1}}}\inps{\frac{1-(-1)^{\ell_k}}{4}} \\
		&= \sum_{k:\;\ell_k\mbox{ odd}} N_k\cdot(-2)^{\ell_1+...+\ell_k-1}
		\end{align}
		is, if at least one $\ell_k$ is odd, equivalent to $-2$ being the root of a non-zero polynomial with only odd coefficients.  This is impossible, contradiction.
	\end{proof}
\end{theorem}
In what follows, let $c^{-1}(N)$ signify the odd preimage of $N$ under the Collatz function.
\begin{theorem}\label{thm:negoneroot}
	If $N = c^{-1}(2^t)$, where $t$ is the base of $N$, and $N$ is not a power of two, then $P_N(-1) = 0$.
	\begin{proof}
		\begin{align}
		P_N(-1) &= \frac{2^{t+1} - 1}{3} + 2^{t+1}\cdot\sum_{k = 1}^{t+1} 2^{-k}\inps{-1}^{k} \\
		&= \frac{2^{t+1} - 1}{3} + 2^{t+1}\cdot\inps{-\frac{1}{2}}\cdot\frac{1 - \inps{-\frac{1}{2}}^{t+1}}{1 - \inps{-\frac{1}{2}}} \\
		&= \frac{2^{t+1} - 1}{3} - \frac{1}{3}\cdot2^{t+1}\inps{1-(-2)^{-t-1}}
		\intertext{Since the base of $N$ is always odd if $N$ is not a power of 2, $(-2)^{-t-1} = 2^{-t-1}$ and this expression simplifies to}
		&= \frac{2^{t+1} - 1}{3} - \frac{2^{t+1}-1}{3} = 0
		\end{align}
	\end{proof}
\end{theorem}
\begin{theorem}
	If $$P_N(-1) = 0$$ then $m(N)$ is even.
	
	\begin{proof}
		If $$P_N(-1) = 0$$ then $$P_N(1) = 0$$ in $\inps{\Z / 2\Z}[z]$, but this is true if and only if $P_N(z)$ has an even number of odd coefficients.
	\end{proof}
\end{theorem}
Unfortunately, a full converse of Theorem \ref{thm:negoneroot} has proved elusive.  For example, 
$$P_{820569}(-1) = 0$$
yet $c(820569) = 1230854 = 2\cdot 615427$, which is not a power of two; in fact, a consequence of Lemma \ref{lemma:fourn} is that
$$P_N(-1) = 0$$ for every odd preimage $N$ of $2^k\cdot 615427$ where $k$ is odd.
\begin{lemma} \label{lemma:fourn}
	If $P_{c^{-1}(N)}(-1) = 0$ then $P_{c^{-1}(4N)}(-1) = 0$.
	\begin{proof}
		If $P_{c^{-1}(N)}(-1) = 0$ then $P_N(-1) = c^{-1}(N) =  \frac{2N-1}{3}$.  But then
		\begin{align}
			P_{c^{-1}(4N)} &= c^{-1}(4N) - 4N + 2N - P_{N}(-1) \\
			&= \frac{8N - 1}{3} - \frac{12N}{3} + \frac{6N}{3} - \frac{2N - 1}{3} \\
			&= 0
		\end{align}
	\end{proof}
\end{lemma}
Finally, there exist $N$ even such that $P_N(-1) = 0$.  The least such example is 
\begin{align}
	N &= \phantom{0}6094358 = 2\cdot 83 \cdot 36713
\intertext{Another example is}
	N &= 46507804 = 2^2 \cdot 7 \cdot 593 \cdot 2801
\end{align}
These two prime factorizations are apparently unrelated.

\section{Conclusion}
Having proved several general properties of Collatz polynomial zeros, including upper and lower bounds, we hope to continue our research along several avenues.  First, we would like to improve the determinantal bounds proved here, perhaps by finding a way of calculating the quantity $U_m$ of Equation \eqref{eq:polbnd} for general large $m$.  Second, we seek to uncover further connections, similar to the conclusions drawn in analytic combinatorics for meromorphic generating functions generally, between the nature of the zeros of Collatz polynomials and the behavior of the Collatz dynamical system.  Third, we would like to consider applications of the unique (and likely difficult to calculate) factorization of the integer $N$ that the $N^{th}$ Collatz polynomial provides by its unique decomposition into linear factors in $\C[z]$; perhaps, for example, this unique factorization could serve as the basis for new cryptographic methods.  Finally, we mention that Collatz polynomials provide a good source for polynomial root-finding algorithms, especially those that seek to find all the roots of a polynomial.

\bibliography{Collatz-Final-arxiv}

\end{document}